\documentclass[10pt,twoside]{article}
\usepackage{amssymb}
\usepackage{latexsym}
\newtheorem{theorem}{Theorem}[section]
\newtheorem{definition}{Definition}[section]
\newtheorem{lemma}{Lemma}[section]

\makeatletter\def\theequation{\arabic{section}.\arabic{equation}}\makeatother
\def\qed{\hfill $\square$}
\def\squar{\vbox{\hrule\hbox{\vrule height 6pt \hskip
6pt\vrule}\hrule}}

\def\qed{\hfill $\squar$}
\def\squar{\vbox{\hrule\hbox{\vrule height 6pt \hskip
6pt\vrule}\hrule}}
\newcommand{\eps}{\varepsilon}

\begin{document}
\title{
\vspace{0.0in} {\bf\Large  Periodic bifurcation from families of
periodic solutions for semilinear differential equations with
Lipschitzian perturbations in Banach spaces}}
\author{{\bf\large Mikhail Kamenskii}\footnote{The author acknowledges the
support by RFBR Grants 05-01-00100,
06-01-72552.}\hspace{2mm}\\
{\it\small Department of Mathematics},\\ {\it\small Voronezh State University}, {\it\small 394006 Voronezh, Russia}\\
{\it\small e-mail: mikhailkamenski@mail.ru}\vspace{1mm}\\
{\bf\large Oleg Makarenkov}\footnote{The author acknowledges the
support by the Grant BF6M10 of Russian Federation Ministry of
Education and CRDF (US), and by RFBR Grants 07-01-00035,
06-01-72552.}\hspace{2mm}\\
{\it\small Research Institute of Mathematics},\\ {\it\small Voronezh State University}, {\it\small 394006 Voronezh, Russia}\\
{\it\small e-mail: omakarenkov@math.vsu.ru}\vspace{1mm}\\
{\bf\large Paolo Nistri}\footnote{The author acknowledges the
support by  INdAM-GNAMPA and the National Research Project PRIN
``Mathematical Control Theory: Controllability,
Optimization, Stability''.}\hspace{2mm}\\
{\it\small Dipartimento di Ingegneria dell'Informazione},\\
{\it\small Universit\`a di Siena}, {\it\small 53100 Siena, Italy}\\
{\it\small e-mail: pnistri@dii.unisi.it}}
\date{}

\maketitle
\begin{center}
{\bf\small Abstract.}

\vspace{2mm} \hspace{.05in}\parbox{4.5in} {{\small Let $A:D(A)\to
E,$ $D(A)\subset E,$ be an infinitesimal generator either of an
analytic compact semigroup or of a contractive $C_0$-semigroup of
linear operators acting in a Banach space $E.$ In this paper we
give both necessary and sufficient conditions for bifurcation of
$T$-periodic solutions for the equation $\dot x=Ax+f(t,x)+\eps
g(t,x,\eps)$ from a $k$-parameterized family of $T$-periodic
solutions of the unperturbed equation corresponding to $\eps=0$.
We show that by means of a suitable modification of the classical
Mel'nikov approach we can construct a bifurcation function and to
formulate the conditions for the existence of bifurcation in terms
of the topological index of the bifurcation function. To do this,
since the perturbation term $g$ is only Lipschitzian we need to
extend the classical Lyapunov-Schmidt reduction to the present
nonsmooth case.}}
\end{center}


\section{\bf Introduction}
\def\theequation{1.\arabic{equation}}\makeatother
\setcounter{equation}{0}
The aim of this paper is to give both
necessary and sufficient conditions for the bifurcation of
$T$-periodic solutions of the semi-linear differential equation
\begin{equation}\label{ps}
  \dot x=Ax+f(t,x)+\eps g(t,x,\eps)
\end{equation}
from a $k$-parameterized family of $T$-periodic solutions of the
unperturbed system, obtained from (\ref{ps}) by letting $\eps=0.$
Here $A:D(A)\to E,$ $D(A)\subset E,$ is an infinitesimal generator
either of an analytic compact semigroup or of a contractive
$C_0$-semigroup of linear operators acting in the Banach space $E,$
the nonlinear operators
 $f\in C^1(\mathbb{R}\times E,E)$ and $g\in
C^0(\mathbb{R}\times E\times[0,1],E)$ are $T$-periodic in the first
variable.

In the case when the unperturbed system is autonomous the problem
was studied by Henry in (\cite{hen}, Ch.~8), where it is assumed
that $g$ is differentiable in the second variable. The author
provided sufficient conditions for bifurcation of $T$-periodic
solutions from a $T$-periodic cycle $x_0,$ the main tool employed in
that paper is the classical Lyapunov-Schmidt reduction, see for
instance Chow and Hale (\cite{ch}, Ch.~2, \S~4). These conditions
are formulated in terms of the existence of nondegenerate zeros of
an analogue of the Malkin's bifurcation function \cite{mal} for an
infinite dimensional Banach space.

In the finite dimensional case, using topological degree
arguments, Felmer and Man\'asevich in \cite{fel} replaced the
assumption of the existence of nondegenerate  zeros of the
bifurcation function by the request that the topological degree of
the bifurcation function is different from zero with respect to a
suitable set. Starting from \cite{fel} there has been a great
amount of work for developing bifurcation results by using the
topological degree theory, see e.g. Henrard and Zanolin \cite{zan}
for bifurcation from a cycle of a Hamiltonian system and
Kamenskii, Makarenkov and Nistri \cite{nach} for bifurcation from
a cycle of a self-oscillating system.  In the present paper we
avoid the requirement that the zeros of the bifurcation function
are nondegenerate, instead we formulate suitable assumptions on
the bifurcation function in terms of the topological degree to
obtain for (\ref{ps}) results similar to those of (\cite{hen},
Ch.~8). To this end we prove an extension of the classical
Lyapunov-Schmidt reduction as presented in (\cite{ch}, Ch.~2,
\S~4) to the case when the perturbation $g$ is Lipschitzian.  We
mention in the sequel some problems involving partial differential
equations which reduce to the situation considered in this paper.
In Chow and Hale \cite[Ch.~8, \S~6]{ch} and Schaeffer and
Golubitsky \cite{gol} the problem of the dependance of the steady
states in chemical reaction models on the relative diffusion
coefficients leads to the consideration of perturbed equations in
Banach spaces with the property that the corresponding unperturbed
equations have a family of solutions.

Another example of such a situation is presented in  Berti and Bolle
\cite{ber}, where the problem of finding periodic solutions of a
nonlinear wave equation by variational methods gives rise to an
unperturbed equation with a family of periodic solutions.

The paper is organized as follows. A modified Lyapunov-Schmidt
reduction for Lipschitzian perturbations  of an operator of the
form $(P-I),$ with $P\in C^1(E,E),$  is obtained in Section~2. In
order to apply the results of Section~2 some relevant properties
of the Poincar\' e map for system (\ref{ps}) are established in
Section~3. Both necessary and sufficient conditions for
bifurcation of periodic solutions to (\ref{ps}) are obtained in
Section~4. Finally, in the appendix of Section~5 we give a proof
of a technical result needed in Section~3.

\section{Lyapunov-Schmidt  reduction}\label{LS}
\def\theequation{2.\arabic{equation}}\makeatother
\setcounter{equation}{0}

\noindent Let $E$ be a Banach space and consider the function
$F:E\times[0,1]\to E$ given by
$$
  F(\xi,\eps)=P(\xi)-\xi+\eps Q(\xi,\eps),
$$
where $P:E\to E$ and $Q:E\times[0,1]\to E.$  Assume that

\begin{itemize}
\item[$(A_1)$] there exist $h_0\in\mathbb{R}^k,$ $r_0>0$
and a function $S\in C^1(B_{\mathbb{R}^k}(h_0,r_0),E)$  such that
$$
  P(\xi)=\xi\quad{\rm for\ any}\quad \xi\in \mathcal{Z}=\{S(h): h\in
  B_{\mathbb{R}^k}(h_0,r_0)\}.
$$
\end{itemize}

\noindent Here and in what follows $B_{X}(c,r)$ denotes the ball in
the normed space $X$ centered at $c$ with radius $r>0$. It is well
known that, under the assumption ($A_1$) with $P\in C^1(E,E)$ and
$Q\in C^1(E\times[0,1],E),$ the Lyapunov-Schmidt reduction
(\cite{ch}, Ch.~2, \S~4) allows to solve the equation
\begin{equation}\label{a}
  F(\xi,\eps)=0,
\end{equation}
for $\eps>0$ sufficiently small.  Next theorem extends this result
to the case when $Q$ satisfies the following Lipschitz condition:

\begin{itemize}
\item[$(L)$]  For any $R>0$ there exists $L(R)>0$ such that
$$
\left\|Q(\xi_1, \eps) - Q(\xi_2, \eps)\right\|\le L(R)
\left\|\xi_1-\xi_2\right\|
$$
whenever $\xi_1, \xi_2\in B_E(0,R)$
and $\eps\in [0,1].$
\end{itemize}

\begin{theorem} \label{th1} Let  $P\in C^1(E,E),$ $Q\in
C^0(E\times[0,1],E),$ where $E$ is a Banach space. Assume that $Q$
satisfies $(L)$. Moreover, assume ($A_1$) and

\begin{itemize}
\item[$(A_2)$] ${\rm dim}S'(h_0)\mathbb{R}^k=k.$
\end{itemize}

\noindent Let $E_{1,h}=S'(h)\mathbb{R}^k.$ Let $E_{2,h}$ be any
subspace of $E$ such that $E=E_{1,h}\bigoplus E_{2,h}$ and assume
that

\begin{itemize}
\item[$(A_3)$] there exists $r_0>0$ such that
both the projectors $\pi_{1,h}$ of $E$ onto $E_{1,h}$ along
$E_{2,h}$ and $\pi_{2,h}$ of $E$ onto $E_{2,h}$ along $E_{1,h}$ are
continuous in $h\in B_{\mathbb{R}^k}(h_0,r_0),$
\end{itemize}

\begin{itemize}
\item[$(A_4)$] for $\xi_0=S(h_0)$ we have
\begin{equation}\label{th1b}
\pi_{2,h_0}(P'(\xi_0)-I)\pi_{2,h_0} \ {\it is\ invertible\ on\ }
E_{2,h_0}.
\end{equation}
\end{itemize}

\noindent Then there exist $0<r_2<r_1<r_0$ and functions
$H:B_E(\xi_0,r_1)\to\mathbb{R}^k,$ with $H(\xi)\to h_0$ as
$\xi\to\xi_0$ and $\beta:B_{\mathbb{R}^k}(h_0,r_1)\times[0,r_1]\to
E,$ $\beta(\cdot,\eps)\in C^0(B_{\mathbb{R}^k}(h_0,r_1),E),$
$\|\beta(h,\eps)\|\le M\eps$ for some $M>0,$ any $h\in
B_{\mathbb{R}^k}(h_0,r_1)$ and any $\eps\in[0,r_1]$ with
\begin{equation}\label{th1a}
  \beta(h,\eps)\in E_{2,h},
\end{equation}
\begin{equation}\label{th1c}
\begin{array}{l}
{\it and\ } {\beta(h,\eps)}/{\eps}\to
-\left(\pi_{2,h}(P'(S(h))-I)\pi_{2,h}\right)^{-1}\pi_{2,h}Q(S(h),0)\quad{\it
as}\ \eps\to 0,\\
{\it uniformly\ in\ }h\in B_{\mathbb{R}^k}(h_0,r_1) \end{array}
\end{equation}
such that the following properties hold:

\vskip0.2truecm 1) if $(\xi,\eps)\in B_E(\xi_0,r_2)\times[0,r_2]$
is a solution to
 equation (\ref{a}) then $(h,\eps),$ where $h=H(\xi),$ is a solution to
\begin{equation}\label{b}
\begin{array}{l}
  \left(S'(h)\right)^{-1}\pi_{1,h}\left[P(\beta(h,\eps)+S(h))\right.\\
  \left.
  -(\beta(h,\eps)+S(h))+\eps
  Q(\beta(h,\eps)+S(h),\eps)\right]=0.
\end{array}
\end{equation}

2) if $(h,\eps)\in B_{\mathbb{R}^k}(h_0,r_1)\times[0,r_1]$ solves
(\ref{b}) then $(\xi,\eps)$ solves (\ref{a}), with
\begin{equation}\label{sol}
  \xi=\beta(h,\eps)+S(h)
\end{equation}

\end{theorem}

Note, that the existence of $\left(S'(h)\right)^{-1}$ on $E_{1,h}$
for $h\in\mathbb{R}^k$ sufficiently close to $h_0$ is guaranteed by
($A_2$) and ($A_3$). To prove Theorem~\ref{th1} we need the
following version of the implicit function theorem.

\begin{lemma}\label{lemimp} Let $E$ be a Banach space and ${V}\subset\mathbb{R}^k$
be an open bounded set. Consider a family of projectors
$\{\pi_h\}_{h\in V}$ on $E$ continuous in $h$ and let $E_h=\pi_h E$
for any $h\in V.$ Let $\Phi_{h,\eps}:E_h\to E_h$ be defined by
\begin{equation}\label{tre}
  \Phi_{h,\eps}(\beta)=\widetilde{P}(h,\beta)+\eps
  \widetilde{Q}(h,\beta,\eps),
\end{equation}
where $\widetilde{P}\in C^0(\mathbb{R}^k\times E,E),$
$\widetilde{Q}\in C^0(\mathbb{R}^k\times E\times[0,1],E),$
$\widetilde{P}(h,\cdot),\widetilde{Q}(h,\cdot,\eps):E_h\to E_h$ for
any $h\in V,$ $\eps\in[0,1].$ Assume that
\begin{enumerate}
\item\label{DST3} the continuity of $\widetilde{P}$ in the first variable is uniform on any bounded
subset of $V\times E,$
\item\label{DST1} $\widetilde{P}$ is differentiable
with respect to the second variable and the derivative is
continuous in $\overline{V}\times E,$
\item\label{DST2} $\widetilde{Q}$ is Lipschitzian in the second variable uniformly
on any bounded subset of $V\times E\times[0,1].$
\item\label{cond1} $\widetilde{P}(h,0)=0$ for any $h\in V,$
\item\label{cond2}
$\pi_h\widetilde{P}'_\beta(h,0):E_h\to E_h$ is invertible for any
$h\in V$ and $(\pi_h\widetilde{P}'_\beta(h,0))^{-1}\pi_h$ is
continuous in $h\in V.$
\end{enumerate}

\noindent Then there exist $r>0,$ $M>0$ and a function
$\beta:V\times[0,r]\to E,$ $\beta(\cdot,\eps)\in C^0(V,E)$
such that

\vskip0.2truecm a) $\beta(h,\eps)\in E_h$ for any $h\in V,$
$\eps\in[0,r],$

\vskip0.2truecm b) $\Phi_{h,\eps}(\beta(h,\eps))=0$ for any $h\in
V,$ $\eps\in[0,r],$

\vskip0.2truecm c) $\beta(h,\eps)$ is the only zero of
$\Phi_{h,\eps}$ in $B_{E_h}(0,r)$ for any $h\in V,$ $\eps\in[0,r],$

\vskip0.2truecm d) $\left\|\beta(h,\eps)\right\|\le M\eps$ for any
$h\in V,$ $\eps\in [0,r].$

\end{lemma}

Although Lemma~\ref{lemimp} looks well-known, the authors were
unable to find a proof of it in the literature, thus for the reader
convenience we provide a proof of Lemma~\ref{lemimp} in the Appendix
of Section~5.

\vskip0.2truecm \noindent {\bf Proof of Theorem~\ref{th1}. }  In
order to define the function $\beta$ we consider the following
auxiliary function $\Phi_{h,\eps}\in C^0(E_{2,h},E_{2,h})$ given by
$$
  \Phi_{h,\eps}(\beta)=\pi_{2,h}\left[P(\pi_{2,h}\beta+S(h))-(\pi_{2,h}\beta+S(h))+\eps
  Q(\beta+S(h),\eps)\right].
$$
Since $P\in C^1(E,E)$ and $S\in C^1(B_{\mathbb{R}^k}(h_0,r_0),E)$
then assumptions \ref{DST3} and \ref{DST1} of Lemma~\ref{lemimp} are
satisfied.

\noindent By our assumptions we have that the application
$(h,\beta,\eps)\to\Phi_{h,\eps}(\beta)$ is Lipschitzian in $\beta$
uniformly on any bounded subset of
$B_{\mathbb{R}^k}(h_0,r_0)\times E\times[0,1]$ and taking into
account ($A_1$) we have

\begin{itemize}
\item[$1)$] $\Phi_{h,0}(0)=0$ for any $h\in B_{\mathbb{R}^k}(h_0,r_0).$
\end{itemize}

\noindent By assumptions $(A_3)$-$(A_4)$ $r_0>0$ can be diminished
in such a way that

\begin{itemize}
\item[$2)$]
$\left(\Phi_{h,0}\right)'(0)=\pi_{2,h}\left(P'(S(h))-I\right)\pi_{2,h}$
is an invertible operator from $E_{2,h}$ to $E_{2,h}$ for $h\in
B_{\mathbb{R}^k}(h_0,r_0).$
\end{itemize}

\noindent Therefore, Lemma~\ref{lemimp} applies with
$$\widetilde{P}(h,\beta)=\pi_{2,h}[P(\pi_{2,h}\beta+S(h))-(\pi_{2,h}\beta+S(h))],$$
$$ \widetilde{Q}(h,\beta,\eps)=\pi_{2,h}Q(\beta+S(h),\eps) \quad
\mbox{and}\quad V=B_{\mathbb{R}^k}(h_0,r_0).
$$
Thus  there exist $r_1\in[0,r_0],$ $M>0$ and a function
$\beta(\cdot,\eps)\in C^0(B_{\mathbb{R}^k}(h_0,r_1),E)$ satisfying
Properties a), b), c) and d) of Lemma~\ref{lemimp}. In particular,
from Property~b) we have
\begin{eqnarray*}
 && \pi_{2,h}[P(\beta(h,\eps)+S(h))-(\beta(h,\eps)+S(h))-\\
  &&-(P(S(h))-S(h))+\eps
  Q(\beta(h,\eps)+S(h),\eps)]=0
\end{eqnarray*}
or equivalently
$$
  \pi_{2,h}\left[(P'(S(h))-I)\pi_{2,h}\beta(h,\eps)+o(\beta(h,\eps))+\eps
  Q(\beta(h,\eps)+S(h),\eps)\right]=0,
$$
for any $h\in B_{\mathbb{R}^k}(h_0,r_1).$

\noindent Therefore
$$
  \beta(h,\eps)=-\left(\pi_{2,h}(P'(S(h))-I)\pi_{2,h}\right)^{-1}\left(\pi_{2,h}o(\beta(h,\eps))+\pi_{2,h}\eps
  Q(\beta(h,\eps)+S(h),\eps)\right).
$$
Due to Property~d) the last equation implies (\ref{th1c}).

\noindent We now proceed to define the function $H$. For this by
($A_2$) we have that $r_1>0$ can be taken sufficiently small such
that $S'(h):\mathbb{R}^k\to E_{1,h}$ is invertible. Thus we can
define the function $\Phi_\xi:\mathbb{R}^k\to\mathbb{R}^k,$
$\xi\in E,$ as follows
$$
  \Phi_\xi(h)=\left(S'(h)\right)^{-1}\pi_{1,h}(\xi-S(h)),\quad
  h\in B_{\mathbb{R}^k}(h_0,r_1).
$$

\noindent We have the following properties for $\Phi_{\xi}.$

\begin{itemize}
\item[$1)$] $\Phi_{\xi_0}$ is differentiable at $h_0.$

\item[$2)$]
$\left(\Phi_{\xi_0}\right)'(h_0)=\left(S'(h_0)\right)^{-1}\pi_{1,h_0}\left(-S'(h_0)\right)=-I,$
namely $\left(\Phi_{\xi_0}\right)'(h_0)$ is an invertible $k\times
k$-matrix.
\end{itemize}

\noindent Observe that property 1) is a direct consequence of the
fact that $\xi_0-S(h_0)=0$ and the continuity of the function
$h\to S^{-1}(h)\pi_h,$ therefore the differentiability of
$\pi_{1,h}$ at $h=h_0$ is not necessary for the validity of 1).

\noindent Let $\delta>0$ be such that $h_0$ is the only zero of
$\Phi_{\xi_0}$ in $B_{\mathbb{R}^k}(h_0,\delta).$ By
(\cite{krazab}, Theorem~6.3) we can consider $\delta>0$
sufficiently small in such a way that
$d(\Phi_{\xi_0},B_{\mathbb{R}^k}(h_0,\delta))=(-1)^k.$ By the
continuity property of the topological degree $r_1>0$ can be
diminished, if necessary, in such a way that
$d(\Phi_\xi,B_{\mathbb{R}^k}(h_0,\delta))=(-1)^k$ for any $\xi\in
B_E(\xi_0,r_1).$ Therefore, for any $\xi\in B_E(\xi_0,r_1)$ there
exists $H(\xi)\in B_{\mathbb{R}^k}(h_0,\delta)$ such that
$\Phi_\xi(H(\xi))=0.$ Let us show that $H(\xi)\to h_0$ as
$\xi\to\xi_0.$ Indeed, arguing by contradiction we would have a
sequence $\{\xi_n\}_{n\in\mathbb{N}}\subset B_E(\xi_0,r_1),$
$h_*\in B_{\mathbb{R}^k}(h_0,\delta)$ such that $H(\xi_n)\to
h_*\not=h_0$ as $n\to\infty$ and thus $\Phi_{\xi_0}(h_*)=0$
contradicting the choice of $\delta>0.$ Therefore
\begin{equation}\label{ro}
\left.\begin{array}{l} \pi_{1,H(\xi)}(\xi-S(H(\xi)))=0,\quad \xi\in
B_E(\xi_0,r_1).
\end{array}\right.
\end{equation}

\noindent Moreover, we consider $r_2\in(0,r_1]$ sufficiently small
to have
\begin{equation}\label{oo}
\left\|\xi-S(H(\xi))\right\|\le r_1,\quad \xi\in B_E(\xi_0,r_2).
\end{equation}

\noindent We are now in the position to complete the proof. For this
let $(\xi,\eps)\in B_E(\xi_0,r_2)\times[0,r_2]$ satisfying
(\ref{a}). Then $(\xi,\eps)$ also satisfies
\begin{eqnarray*}
\left\{\begin{array}{l}
\pi_{1,H(\xi)}[P\left(\xi-S(H(\xi))+S(H(\xi))\right)-\\-\left(\xi-S(H(\xi))+S(H(\xi))\right)+
\eps Q\left(\xi-S(H(\xi))+S(H(\xi)),\eps\right)]=0,\\
\pi_{2,H(\xi)}[P\left(\xi-S(H(\xi))+S(H(\xi))\right)-\\-\left(\xi-S(H(\xi))+S(H(\xi))\right)+
\eps Q\left(\xi-S(H(\xi))+S(H(\xi)),\eps\right)]=0.
\end{array}\right.
\end{eqnarray*}
From (\ref{ro}), (\ref{oo}) and  Property~c) of Lemma~\ref{lemimp}
we have
\begin{eqnarray}\label{zz}
\left\{\begin{array}{l}
\pi_{1,H(\xi)}[P\left(\xi-S(H(\xi))+S(H(\xi))\right)-
\\ -\left(\xi-S(H(\xi))+S(H(\xi))\right)+
\eps Q\left(\xi-S(H(\xi))+S(H(\xi)),\eps\right)]=0,\\
\beta(H(\xi),\eps)=\xi-S(H(\xi)).
\end{array}\right.
\end{eqnarray}
Therefore,
\begin{equation}\label{zzz_}
  \pi_{1,h}\left[P(\beta(h,\eps)+S(h))-(\beta(h,\eps)+S(h))+\eps
  Q(\beta(h,\eps),\eps)\right]=0
\end{equation}
has a solution $h=H(\xi).$ Since $r_1>0$ has been chosen in such a
way that $S'(h)$ is invertible on $E_{1,h}$ for $h\in
B_{\mathbb{R}^k}(h_0,r_1)$ then (\ref{zzz_}) can be rewritten as
(\ref{b}). Assume now that (\ref{b}) is satisfied with some
$(h_*,\eps_*)\in B_{\mathbb{R}^k}(h_0,r_1)\times[0,r_1].$ Define
$\widetilde{\xi}\in E$ as
\begin{equation}\label{ti}
  \widetilde{\xi}=\beta(h_*,\eps_*)+S(h_*).
\end{equation}
Since $(S'(h_*))^{-1}$ is invertible then
$\pi_{1,h_*}[P(\widetilde{\xi})-\widetilde{\xi}+\eps
Q(\widetilde{\xi},\eps_*)]=0.$ On the other hand from (\ref{ti})
we have
$$
\begin{array}{l}
\pi_{2,h_*}[P(\pi_{2,h_*}\beta(h_*,\eps_*)+S(h_*))-(\pi_{2,h_*}\beta(h_*,\eps_*)+S(h_*))+\\
 +\eps Q(\beta(h_*,\eps_*)+S(h_*),\eps)]
 =\pi_{2,h_*}[P(\widetilde{\xi})-\widetilde{\xi}+\eps
Q(\widetilde{\xi},\eps)].
\end{array}
$$
Thus $(\xi_*,\eps_*)$ solves (\ref{a}) and so the proof is
complete.\qed

\vskip0.4truecm\noindent The following two results are consequences
of Theorem~\ref{th1} and they provide, respectively, a necessary and
a sufficient condition for the existence of solutions to (\ref{a})
near $\xi_0$ when $\eps>0$ is sufficiently small. These conditions
are expressed in terms of the following bifurcation function
\begin{eqnarray*}
M(h)&=&\left(S'(h)\right)^{-1}\pi_{1,h}[Q(S(h_0),0)-\\&&-
\left(P'(S(h))-I\right)\left(\pi_{2,h}(P'(S(h))-I)\pi_{2,h}\right)^{-1}\pi_{2,h}Q(S(h),0)],
\end{eqnarray*}
where $h$ varies in a sufficiently small neighborhood of
$h_0\in\mathbb{R}^k.$

\vskip0.2truecm We can prove the following. \vskip0.2truecm
\begin{theorem} \label{thm2} Let all the assumptions of Theorem~\ref{th1}
be satisfied. Assume that there exist sequences $\eps_n\to 0$ and
$\xi_n\to \xi_0$ as $n\to\infty$ such that $(\xi_n,\eps_n)$ solves
(\ref{a}). Then
\begin{equation}\label{nec}
M(h_0)=0.
\end{equation}
\end{theorem}

\noindent {\bf Proof.} By Theorem~\ref{th1}, for $n\ge n_0,$ with
$n_0\in\mathbb{N}$ sufficiently large, we have that
\begin{eqnarray}\label{n1}
 && \left(S'(h_n)\right)^{-1}\pi_{1,h_n}[P(\beta(h_n,\eps_n)+S(h_n))-\nonumber\\
 && -(\beta(h_n,\eps_n)+S(h_n))+\eps_n
  Q(\beta(h_n,\eps_n)+S(h_n),\eps_n)]=0
\end{eqnarray}
where $h_n=H(\xi_n)$. On the other hand $n_0$ can be chosen
sufficiently large in such a way that
$$
  P\left(S(h_n)\right)-S(h_n)=0\quad{\rm for\ } n\ge n_0
$$
thus, for $n\ge n_0,$ (\ref{n1}) can be rewritten as
\begin{eqnarray}\label{sprime}
  &&\left(S'(h_n)\right)^{-1}\pi_{1,h_n}[\left(P'(S(h_n))-I\right)\frac{\beta(h_n,\eps_n)}{\eps_n}+\nonumber\\
  &&+\frac{o(\beta(h_n,\eps_n))}{\eps_n}+
  Q(\beta(h_n,\eps_n)+S(h_n),\eps_n)]=0.
\end{eqnarray}
By means of property (\ref{th1c}) we can pass to the limit as
$n\to\infty$ in (\ref{sprime}) to obtain (\ref{nec}). \qed

\begin{theorem} \label{thm3} Let all the assumptions of Theorem~\ref{th1} be
satisfied. Assume that
\begin{equation}\label{th3a}
  h_0\ {\rm is\ an\ isolated\ zero\ of\ }M
\end{equation}
and
\begin{equation}\label{th3b}
  {\rm ind}\left(h_0,M\right)\not=0.
\end{equation}
Then, for any $\eps>0$ sufficiently small there exists
$\xi_\eps\in E$ such that
$$
  F(\xi_\eps,\eps)=0
$$
and
 \begin{equation}\label{conve}
    \xi_\eps\to\xi_0\quad{\rm as }\ \eps\to 0.
 \end{equation}
\end{theorem}

\noindent {\bf Proof.} Let $r_1>0$ be as given by Theorem~\ref{th1}.
Since
\begin{equation}\label{the}
 P(S(h))=S(h)\quad{\rm for\ any\ }h\in B_{\mathbb{R}^k}(h_0,r_1)
\end{equation}
then the zeros of the function
\begin{eqnarray*}
 \Phi(h,\eps)&=&\left(S'(h)\right)^{-1}\pi_{1,h}[P(\beta(h,\eps)+S(h))-\\&&-(\beta(h,\eps)+S(h))+\eps
  Q(\beta(h,\eps)+S(h),\eps)]
\end{eqnarray*}
coincide with the zeros of the function
\begin{eqnarray*}
 M_\eps(h)&=&\left(S'(h)\right)^{-1}\pi_{1,h}[\left(P'(S(h))-I\right)\frac{\beta(h,\eps)}{\eps}+\\&&+
  \frac{o(\beta(h,\eps))}{\eps}+
  Q(\beta(h,\eps)+S(h),\eps)].
\end{eqnarray*}
In order to apply Theorem~\ref{th1} we show now that $r\in(0,r_1]$
can be chosen in such a way that the function $M_\eps$ has zeros in
$B_{\mathbb{R}^k}(h_0,r)$ for any $\eps>0$ sufficiently small.

\vskip0.2truecm\noindent By condition (\ref{th3a}) $r>0$ can be
chosen sufficiently small in such a way that
\begin{equation}\label{iso1}
 {\rm the\ only\ zero\ of\ }M \ {\rm in\ } B_{\mathbb{R}^k}(h_0,r) {\ \rm is\ }h_0.
\end{equation}
Therefore, by condition (\ref{th3b}) we have
$$
  d(M,B_{\mathbb{R}^k}(h_0,r))={\rm ind}(h_0,M)\not=0.
$$
On the other hand from property (\ref{th1c}) we have that
\begin{equation}\label{conv}
  M_\eps(h)\to M(h)\quad{\rm as\ }\eps\to 0
\end{equation}
uniformly with respect to $h\in B_{\mathbb{R}^k}(h_0,r).$ Thus we
conclude that
$$
  d(M_\eps,B_r(h_0))\not=0
$$
for $\eps\in(0,\eps_0],$ where $\eps_0>0$ is sufficiently small.
Thus for any $\eps\in(0,\eps_0]$ there exists $h_\eps$ such that
$M_\eps(h_\eps)=0.$ Moreover, we have that
$$
  h_\eps\to h_0\quad {\rm as}\ \eps\to 0
$$
otherwise $M$ would have zeros in $B_{\mathbb{R}^k}(h_0,r)$
different from $h_0,$ contradicting (\ref{iso1}). Finally,
(\ref{conve}) follows from (\ref{sol}). \qed

\vskip0.2truecm In finite dimensional spaces results similar to
previous Theorems~\ref{thm2} and \ref{thm3} have been recently
obtained by Buica, Llibre and Makarenkov \cite{blm}, where the
uniqueness of the bifurcating periodic solutions is also proved.

\section{The Poincar\' e map}
\def\theequation{3.\arabic{equation}}\makeatother
\setcounter{equation}{0}

Since the definition of the Poincar\' e map for system (\ref{ps}) on
the time interval $[0, T]$ depends on the assumptions on the linear
unbounded operator $A,$ we precise in $(C1)$ and $(C2)$ below the
two cases that we consider for $A$ in the paper.

\begin{itemize}
\item[\bf $(C1)$] The operator $A$ is a generator of an
analytic compact semigroup ${\rm e}^{At}$ in $E.$ The operators
$f,g$ are subordinated to some $A^{-\alpha},$ $0<\alpha<1$ (see
e.g. \cite{KZPS}), the operator $f(\cdot,A^{-\alpha}\cdot)$ is
differentiable in the second variable and the operators
$f'_{(2)}(\cdot,A^{-\alpha}\cdot)$, $g(\cdot,A^{-\alpha}\cdot,
\cdot)$ are continuous in $\mathbb{R}\times E$ and they satisfy a
Lipschitz condition in the second variable uniformly with respect
to the others.
\item[\bf $(C2)$] The operator $A$ is a generator of a
$C_0$-semigroup ${\rm e}^{At}.$ The semigroup ${\rm e}^{At}$ is
contractive, namely
$$
\left\|{\rm e}^{At}\right\|\le{\rm
e}^{-\gamma t},
$$
where $\gamma>0.$ The operators $f$ and $g$ are continuous from
$\mathbb{R}\times E\to E$ and verify the inequality
$$
  \chi(f(t,\Omega))\le k\chi(\Omega),\ \chi(g(t,\Omega,\eps))\le k\chi(\Omega),
$$
where $\chi$ is the Hausdorff measure of noncompactness
\footnote{We recall (see \cite{AKRPS}) that for a bounded set
$\Omega\subset E$ the Hausdorff measure of noncompactness is
defined by the formula
$$
\chi(\Omega)={\rm inf}\left\{r>0:\; \mbox{there exists}\;
(y_1,...,y_m)\; \mbox{such that} \; \Omega\subset\cup_{i=1}^m
B(y_i,r)\right\},
$$
where $m\in \mathbb{N}.$

The continuous operator $F: E\to E$ is called $(q,\chi)$-condensing
if
$$
\chi(F(t,\Omega))\le q\chi(\Omega)
$$
for any bounded $\Omega\in E.$} in the space $E,$ $k\ge 0$ and
$q=k/\gamma<1.$ The operator $f$ is differentiable in the second
variable and the operators $f'_{(2)}$ and $g$ are continuous in
$\mathbb{R}\times E$ and they satisfy a Lipschitz condition in the
second variable uniformly with respect to the others.
\end{itemize}

\noindent It is a classical result (see e.g. \cite{KZPS}) that
$(C1)$ and $(C2)$ ensures respectively that the integral equations
\begin{equation}\label{i1}
 x(t)={\rm e}^{At}\xi+\int\limits_0^t A^\alpha{\rm
 e}^{A(t-s)}\left[f(s,A^{-\alpha}x(s))+\eps
 g(s,A^{-\alpha}x(s),\eps)\right]ds,
\end{equation}

\begin{equation}\label{i2}
 x(t)={\rm e}^{At}\xi+\int\limits_0^t {\rm
 e}^{A(t-s)}\left[f(s,x(s))+\eps g(s,x(s),\eps)\right]ds
\end{equation}
have a unique solution $x(\cdot)$ defined on some interval $[0,d],
d>0.$ By means of this function $x$ we can define the shift operator
as follows.

\begin{definition}
{\rm Let $x:[0, d]\times E\times[0,1]\to E$ be defined at
$(t,\xi,\eps)$ as $x(t,\xi,\eps)=x(t)$ for all $t\in[0,d].$ If for
some $\xi\in E$ and $\eps\in[0,1]$ we have that
$x(\cdot,\xi,\eps)$ is defined on the whole time interval $[0,T]$
then for these values $\xi$ and $\eps$ we define the Poincar\' e
map for system (\ref{ps}) as
$$
\mathcal{P}_\eps(\xi)=x(T,\xi,\eps).
$$}
\end{definition}

\vskip0.2truecm A crucial role in what follows is played by the
following technical lemma.

\begin{lemma}\label{lem2} Assume that either $(C1)$ or $(C2)$ is satisfied.
Assume that for some $\xi_0\in E$ the shift operator $(t,\xi,\eps)
\to x(t,\xi,\eps)$ is well defined for $t=T,$ $\xi=\xi_0$ and
$\eps=0.$ Then there exists $r>0$ such that this operator is well
defined for $t=T,$ any  $\xi\in B_E(\xi_0,r),$ any $\eps\in[0,r]$
and the function
$$
  u(t,\xi,\eps)=\frac{x(t,\xi,\eps)-x(t,\xi,0)}{\eps}
$$
is Lipschitz in the second variable uniformly in $[0,T]\times
B_E(\xi_0,r)\times(0,r],$ namely there exists $L>0$ such that
$$
  \left\|u(t,\xi_1,\eps)-u(t,\xi_2,\eps)\right\|\le
  L\,\|\xi_1-\xi_2\|
$$
for any $t\in [0,T],\ \xi_1,\xi_2\in
  B_E(\xi_0,r)$  and  $\eps\in(0,r].$
\end{lemma}

\noindent {\bf Proof.} The fact that the assumptions of the Lemma
imply the existence of $r>0$ such that the operator
$(t,\xi,\eps)\mapsto x(t,\xi,\eps)$ is well defined, bounded and
continuous on $[0,T]\times B_E(\xi_0,r)\times[0,r]$ is well known,
see, for instance, (\cite{koz}, Theorem~5.2.5). In the sequel we
have $\alpha\not=0$ if $(C1)$ holds, while $\alpha=0$ if we assume
$(C2)$. Since $A^{-\alpha}$ is either a compact operator or the
identity then the operator $(t,\xi,\eps)\mapsto
A^{-\alpha}x(t,\xi,\eps)$ is well defined, bounded and continuous on
$[0,T]\times B_E(\xi_0,r)\times[0,r]$. Therefore, taking into
account that $f'_x$ satisfies Lipschitz condition, there exists
$\widehat{M}>0$ such that
$$\|f'_x(s,A^{-\alpha}\{\theta
x(s,\xi_1,\eps)+(1-\theta)x(s,\xi_2,\eps)\})\|\le\widehat{M}$$ for
any $s\in[0,T],$ $\theta\in[0,1],$ $\xi_1,\xi_2\in B_E(\xi_0,r).$

\vskip0.2truecm\noindent From the continuous differentiability of
$f$ and the Lipschitz condition on $g$ assumed in $(C1)$ and $(C2)$
we deduce the existence of $\widetilde{M}>0$ such that
$$\|f(t,A^{-\alpha}\xi)\|+\|g(t,A^{-\alpha}\xi,\eps)\|\le
\widetilde{M}$$ for any $t\in[0,T],$ $\xi\in
x([0,T],B_E(\xi_0,r),[0,r])$ and $\eps\in[0,r].$

\vskip0.2truecm\noindent Since
$A^{-\alpha}x([0,T],B_E(\xi_0,r),[0,r])$ is bounded then by using
the Lipschitz condition on $g$ we obtain the existence of
$\widehat{L}>0$ such that
$$
\|g(s,A^{-\alpha}\xi_1,\eps)-g(s,A^{-\alpha}\xi_2,\eps)\|\le
\,\widetilde{L}\|\xi_1-\xi_2\|
$$
for any $s\in[0,T],$ $\xi_1,\xi_2\in
x([0,T],B_E(\xi_0,r),[0,r])$ and $\eps\in[0,r].$

\vskip0.2truecm\noindent Furthermore, by \cite[Theorem~6.13]{Pazy}
there exists $c>0$ such that $\sup\limits_{t\in[0,T]}\left\|{\rm
e}^{At}\right\|<c$ and $\left\|A^\alpha{\rm
e}^{At}\right\|<{c}/{t^\alpha},$ where either $\alpha=0$ or
$\alpha>0.$

\vskip0.2truecm\noindent Now given an arbitrary $\phi\in
B_{E^*}(0,1),$ where $E^*$ denotes the dual space of $E$, we
evaluate $\left<\phi,x(t,\xi_1,\eps)-x(t,\xi_2,\eps)\right>$ as
follows

\begin{eqnarray}\label{I}
  &&\left<\phi,x(t,\xi_1,\eps)-x(t,\xi_2,\eps)\right>=\nonumber\\&=&\left<\phi,{\rm
  e}^{At}(\xi_1-\xi_2)\right>+ \int\limits_0^t \big<\phi,A^\alpha {\rm
  e}^{A(t-s)}f'_x\Big(s,A^{-\alpha}\{\theta(s,\xi_1,\xi_2,\eps)x(s,\xi_1,\eps)+\nonumber\\
  &&+(1-\theta(s,\xi_1,\xi_2,\eps)x(s,\xi_2,\eps))\}\Big)\,
  A^{-\alpha}\big(x(s,\xi_1,\eps)-x(s,\xi_2,\eps)\big)\big>ds+\nonumber\\
  &&+\eps\int\limits_0^t\left<\phi,A^\alpha{\rm
  e}^{A(t-s)}\left(g\left(s,A^{-\alpha}
  x(s,\xi_2,\eps),\eps\right)-g\left(s,A^{-\alpha}
  x(s,\xi_1,\eps),\eps\right)\right)\right>ds\le\nonumber\\
  &\le&   c\|\xi_1-\xi_2\|+\int\limits_0^t\frac{c\widehat{M}}{(t-s)^\alpha}\left\|x(s,\xi_1,\eps)-x(s,\xi_2,\eps)\right\|ds
  +\nonumber\\&&+ \eps\int\limits_0^t\frac{c\widetilde{L}}{(t-s)^\alpha}\left\|x(s,\xi_1,\eps)-x(s,\xi_2,\eps)\right\|ds.
\end{eqnarray}
Since $\phi$ is arbitrary we have
\begin{eqnarray}\label{2SSS}
  \left\|x(t,\xi_1,\eps)-x(t,\xi_2,\eps)\right\|
 &\le&
 {c\|\xi_1-\xi_2\|}+\int\limits_0^t\frac{c\widehat{M}}{(t-s)^\alpha}\left\|x(s,\xi_1,\eps)-x(s,\xi_2,\eps)\right\|ds+\nonumber\\&&+
 \eps\int\limits_0^t\frac{c\widetilde{L}}{(t-s)^\alpha}\left\|x(s,\xi_1,\eps)-x(s,\xi_2,\eps)\right\|ds.
\end{eqnarray}
Dividing the last inequality by $\|\xi_1-\xi_2\|$ one obtains that
\begin{eqnarray*}
  \frac{\left\|x(t,\xi_1,\eps)-x(t,\xi_2,\eps)\right\|}{\|\xi_1-\xi_2\|}\le c+\int\limits_0^t\frac{c\widehat{M}+\eps c
  \widetilde{L}}{(t-s)^\alpha}\cdot \frac{\left\|x(s,\xi_1,\eps)-x(s,\xi_2,\eps)\right\|}{\|\xi_1-\xi_2\|}ds.
\end{eqnarray*}
Using the generalized Gronwall--Bellman lemma, see (\cite{hen},
Lemma~7.1.1), from the last inequality we obtain that there exists
$M_v>0$ such that
\begin{equation}\label{M_v}
\frac{\left\|x(t,\xi_1,\eps)-x(t,\xi_2,\eps)\right\|}{\|\xi_1-\xi_2\|}
\le M_v
\end{equation}
for any $(t,\xi_1,\xi_2,\eps)\in [0,T]\times B_E(\xi_0,r)\times
B_E(\xi_0,r)\times[0,r].$

\vskip0.2truecm\noindent For the function $u(t,\xi,\eps)$ we have
the following inequality
\begin {eqnarray*}
&&\left<\phi,u(t,\xi,\eps)\right>=\\
  &=&\Big<\phi,\frac{1}{\eps}\int_0^t A^\alpha {\rm
  e}^{A(t-s)}\left[f(s,A^{-\alpha}x(s,\xi,\eps))-f(s,A^{-\alpha}x(s,\xi,0))\right] ds+\\
 &&\;\;+\int_0^t A^\alpha{\rm
  e}^{A(t-s)}g(s,A^{-\alpha}x(s,\xi,\eps))ds \Big>\le\\
  &\le&
  \int_0^t\frac{c\widehat{M}\|u(s,\xi,\eps)\|}{(t-s)^\alpha}ds+\int_0^t\frac{c\widetilde{M}}{(t-s)^\alpha}ds.
\end{eqnarray*}
Using again the generalized Gronwall--Bellman lemma from the last
inequality we obtain that there exists $M_u>0$ such that
\begin{equation}\label{M_u} \|u(t,\xi,\eps)\|\le
M_u\quad{\rm for\ any\ } (t,\xi,\eps)\in [0,T]\times
B_E(\xi_0,r)\times[0,r].
\end{equation}

\noindent Observe that if a function $\Psi:E\to E$ is
differentiable and there exists $L>0$ such that
$\|\Psi'(\xi)-\Psi'(\zeta)\|\le L\,\|\xi-\zeta\|$ for any
$\xi,\zeta\in E$ then
\begin {eqnarray}\label{M}
 && \|\Psi(\xi_2)-\Psi(\xi_1)-\Psi(\zeta_2)+\Psi(\zeta_1)\|\le\nonumber\\
 & \le & \sup_{0\le\theta\le
  1}\|\Psi'(\zeta_2+\theta(\xi_2-\zeta_2))\|\,\|\xi_2-\xi_1-\zeta_2+\zeta_1\|+\nonumber\\
 && + L\max\{\|\xi_2-\xi_1\|,\|\zeta_2-\zeta_1\|\}\;\|\xi_1-\zeta_1\|.
\end{eqnarray}
To prove this it is sufficient to consider the real function
$\gamma:[0,1]\to \mathbb{R}$ given by
$$
  \gamma(\tau)=\left<\phi,\Psi(\zeta_2+\tau(\xi_2-\zeta_2)-\Psi(\zeta_1+\tau(\xi_1-\zeta_1))\right>,\quad\tau\in[0,1].
$$
By Lagrange theorem there exists $\theta\in[0,1]$ such that
$$
  \gamma(1)-\gamma(0)=\gamma'(\theta)
$$
and then
\begin{eqnarray*}
 &&
 |\left<\phi,\Psi(\xi_2)-\Psi(\xi_1)-\Psi(\zeta_2)+\Psi(\zeta_1)\right>|=\gamma(1)-\gamma(0)=\gamma'(\theta)=\\
 &=&|\left<\phi,\Psi'(\zeta_2+\theta(\xi_2-\zeta_2))(\xi_2-\zeta_2)-\Psi'(\zeta_1+\theta(\xi_1-\zeta_1))(\xi_1-\zeta_1)\right>|\le\\
&\le &|\left<\phi,\Psi'(\zeta_2+\theta(\xi_2-\zeta_2))(\xi_2-\xi_1-\zeta_2+\zeta_1)\right>|+\\
&&+|\left<\phi,(\Psi'(\zeta_2+\theta(\xi_2-\zeta_2))-\Psi'(\zeta_1+\theta(\xi_1-\zeta_1)))(\xi_1-\zeta_1)\right>|\le\\
&\le &\|\Psi'(\zeta_2+\theta(\xi_2-\zeta_2))\|\,\|\xi_2-\xi_1-\zeta_2+\zeta_1\|+\\
&&+\|\Psi'(\zeta_2+\theta(\xi_2-\zeta_2))-\Psi'(\zeta_1+\theta(\xi_1-\zeta_1))\|\,\|\xi_1-\zeta_1\|\le\\
&\le & \sup_{\theta\in [0,1]}\|\Psi'(\zeta_2+\theta(\xi_2-\zeta_2)\|\,\|\xi_2-\xi_1-\zeta_2+\zeta_1\|+\\
&&  +L\|(1-\theta)\zeta_2+\theta \xi_2-(1-\theta)\zeta_1-\theta
  \xi_1\|\cdot\|\xi_1-\xi_2\|=\\
& = & \sup_{\theta\in
[0,1]}\|\Psi'(\zeta_2+\theta(\xi_2-\zeta_2)\|\,
  \|\xi_2-\xi_1-\zeta_2+\zeta_1\|+\\
 && +L\|(1-\theta)(\zeta_2-\zeta_1)+\theta(\xi_2-\xi_1)\|\,\|\xi_1-\xi_2\|\le\\
 & \le & \sup_{\theta\in
 [0,1]}\|\Psi'(\zeta_2+\theta(\xi_2-\zeta_2)\|\,
  \|\xi_2-\xi_1-\zeta_2+\zeta_1\|+\\
&& +L\max\{\|\xi_2-\xi_1\|,\|\zeta_2-\zeta_1\|\}\|\xi_1-\zeta_1\|.
\end{eqnarray*}

\noindent By the Lipschitz assumption on $f'_x$ there exists
$\widehat{L}>0$ such that
$$\|f'_x(s,A^{-1}\xi_1)-f'_x(s,A^{-1}\xi_2)\|\le\widehat{L}\,\|\xi_1-\xi_2\|$$
for any $s\in[0,T],$ $\xi_1,\xi_2\in x([0,T],B_E(\xi_0,r),[0,r]).$

\noindent Consider now
\begin{eqnarray*}
&&  \frac{u(t,\xi_1,\eps)-u(t,\xi_2,\eps)}{\|\xi_1-\xi_2\|}=
\frac{x(t,\xi_1,\eps)-x(t,\xi_1,0)-x(t,\xi_2,\eps)+x(t,\xi_2,0)}{\eps\|\xi_1-\xi_2\|}=\\
& = & \frac{1}{\eps\|\xi_1-\xi_2\|}\int_0^t A^\alpha {\rm
  e}^{A(t-s)}(f(s,A^{-\alpha}x(s,\xi_1,\eps))-f(s,A^{-\alpha}x(s,\xi_1,0))-\\
&&-f(s,A^{-\alpha}x(s,\xi_2,\eps))+f(s,A^{-\alpha}x(s,\xi_2,0)))ds+\\
&&  +\frac{1}{\|\xi_1-\xi_2\|}\int_0^t A^\alpha {\rm
  e}^{A(t-s)}(g(s,A^{-\alpha}x(s,\xi_1,\eps))-g(s,A^{-\alpha}x(s,\xi_2,\eps)))ds\le\\
& \le &
\sup_{s\in[0,T],\theta\in[0,1]}\|f'_x(s,A^{-\alpha}(x(s,\xi_2,0)+
\theta(x(s,\xi_1,0)-x(s,\xi_2,0))))A^{-\alpha}\|\le\\
& \le & \int_0^t \frac{c}{(t-s)^\alpha}\,
\frac{\|x(s,\xi_1,0)-x(s,\xi_1,\eps)-x(s,\xi_2,0)+
x(s,\xi_2,\eps)\|}{\eps\|\xi_1-\xi_2\|}ds +\\
& &
+\sup_{s\in[0,T]}\max\left\{\frac{\|x(s,\xi_1,0)-x(s,\xi_1,\eps)\|}{\eps},
\frac{\|x(s,\xi_2,0)-x(s,\xi_2,\eps)\|}{\eps}\right\}\cdot\\
& &\cdot
\int_0^t\frac{c\widehat{L}}{(t-s)^{\alpha}}\,\frac{\|x(s,\xi_1,\eps)-x(s,\xi_2,\eps)\|}{\|\xi_1-\xi_2\|}
ds
  +\\ &&+ \int_0^t\frac{c\widetilde{L}}{(t-s)^\alpha}
  \,\frac{\|x(s,\xi_1,\eps)-x(s,\xi_2,\eps)\|}{\|\xi_1-\xi_2\|}ds.
\end{eqnarray*}
By (\ref{M_u}) and (\ref{M_v}) there exists $M>0$ such that the
last inequality can be rewritten as
$$
  \frac{u(t,\xi_1,\eps)-u(t,\xi_2,\eps)}{\|\xi_1-\xi_2\|}
  \le  \int_0^t \frac{c\widehat{M}}{(t-s)^\alpha}\,
\frac{\|u(s,\xi_1,\eps)-u(s,\xi_2,\eps)\|}{\|\xi_1-\xi_2\|}ds+M
$$
and the assertion follows from the generalized Gronwall--Bellman
lemma, see (\cite{hen}, Lemma~7.1.1). \qed

\section{Existence of periodic solutions}
\def\theequation{4.\arabic{equation}}\makeatother
\setcounter{equation}{0}

\noindent In this section we assume that either $(C1)$ or $(C2)$ is
satisfied, moreover we assume the following condition:

\begin{itemize}
\item[$(\widetilde{A}_0)$] the solution $x$ of (\ref{ps})
with $\eps=0$ satisfying $x(0)=\xi_0$ is defined on $[0,T],$ namely
the Poincar\' e map $\mathcal{P}_0$ is defined at $\xi_0.$
\end{itemize}

\noindent Therefore, from Lemma~\ref{lem2} we have that there exists
$r>0$ such that  the Poincar\' e map $\mathcal{P}_\eps$ for system
(\ref{ps}) is defined on $B_E(\xi_0,r)$ for any $\eps\in[0,r]$ and
it has the form
$$
  \mathcal{P}_\eps(\xi)=\mathcal{P}_0(\xi)+\eps Q(\xi,\eps),
$$
where $\mathcal{P}_0$ is differentiable and $Q$ satisfies a
Lipschitz condition in the first variable $\xi$ uniformly on
$B_E(\xi_0,r)\times[0,r].$

\vskip0.2truecm\noindent Letting
$F(\xi,\eps)=\mathcal{P}_\eps(\xi)$ assumptions ($A_1$)-($A_4$) of
Theorem~\ref{th1} can be rewritten as
\begin{itemize}
\item[$(\widetilde{A}_1)$] there exists a function $S\in
C^1({V},E)$ defined on some open neighborhood
${V}\subset\mathbb{R}^k$ of $h_0$  such that $S(h_0)=\xi_0$ and $
  \mathcal{P}_0(\xi)=\xi$ for any $\xi\in \mathcal{Z}=\bigcup\limits_{h\in
  {V}}S(h),$

\item[$(\widetilde{A}_2)$]  ${\rm dim}S'(h_0)\mathbb{R}^k=k.$
\end{itemize}

\noindent Let $E_{1,h}=S'(h)\mathbb{R}^k$ and let $E_{2,h}$ be any
subspace of $E$ such that $E=E_{1,h}\bigoplus E_{2,h}$ and

\begin{itemize}
\item[$(\widetilde{A}_3)$] both the projectors
$\pi_{1,h}$ of $E$ onto $E_{1,h}$ along $E_{2,h}$ and $\pi_{2,h}$ of
$E$ onto $E_{2,h}$ along $E_{1,h}$ are continuous in $h\in V,$

\item[$(\widetilde{A}_4)$] for $\xi_0=S(h_0)$ we have
\begin{equation}\label{th1b}
\pi_{2,h_0}((\mathcal{P}_0)'(\xi_0)-I)\pi_{2,h_0} \ {\rm is\
invertible\ on\ } E_{2,h_0}.
\end{equation}
\end{itemize}
Furthermore, it can be observed that $Q(\xi,0)$ is the value of the
solution of the problem
\begin{equation}\label{pro}
\begin{array}{l}
  \dot y=Ay+f'_x(t,x(t,\xi,0))y+g(t,x(t,\xi,0),0),\\
  y(0)=0
\end{array}
\end{equation}
at time $t=T.$

\noindent To see this, observe that the function $u$ of
Lemma~\ref{lem2} satisfies the following integral equation
\begin{eqnarray*}
  & & \hskip-0.5cm u(t,\xi,\eps)=\int_0^t A^\alpha{\rm e}^{\Lambda
  (t-s)}f'_x(s,A^{-\alpha}x(s,\xi,0))u(s,\xi,\eps)ds+\\
  & & \hskip-0.5cm+\int_0^t A^\alpha {\rm
  e}^{\Lambda(t-s)}\,\frac{o(x(s,\xi,\eps)-x(s,\xi,0))}{\eps}ds+
\int_0^t A^\alpha {\rm
  e}^{\Lambda(t-s)}g(s,A^{-\alpha}x(s,\xi,\eps),\eps)ds
\end{eqnarray*}
and so $u(T,\xi,0)=Q(\xi,0).$  Therefore, we can give an
equivalent definition of the bifurcation function $M$ introduced
in Section~\ref{LS}, that is $M\in C^0(\mathbb{R}^k,\mathbb{R}^k)$
can be defined as follows
\begin{eqnarray*}
 M(h)&=&\left(S'(h)\right)^{-1}\pi_{1,h}[\eta(S(h_0))-\\
&&-((\mathcal{P}_0)'(S(h))-
I)\left(\pi_{2,h}((\mathcal{P}_0)'(S(h))-I)\pi_{2,h}\right)^{-1}\pi_{2,h}\eta(S(h))],
\end{eqnarray*}
where $h\in B_{\mathbb{R}^k}(h_0,r),$ and $\eta$ is the value of the
solution of (\ref{pro}) at time $t=T.$

\vskip0.4truecm From Theorem~\ref{thm2} we have the following
necessary condition for the existence of $T$-periodic solutions to
(\ref{ps}).

\begin{theorem} Assume that $(C1)$ or $(C2)$ is satisfied. Assume
($\widetilde{A}_0$)-($\widetilde{A}_4$). Assume that there exists a
sequence $\eps_n\to 0$ as $n\to\infty$ and a sequence of
$T$-periodic functions $x_n\in C^0([0,T],E),$ $x_n\to
x(\cdot,\xi_0,0)$ as $n\to\infty$ such that $(x_n,\eps_n)$ solves
(\ref{ps}). Then
\begin{eqnarray*}\label{nec1}
M(h_0)=0.
\end{eqnarray*}
\end{theorem}

Analogously from Theorem~\ref{thm3} we derive the following
sufficient condition for the existence of $T$-periodic solutions
to (\ref{ps}).

\begin{theorem}  Assume that $(C1)$ or $(C2)$ is satisfied. Assume
($\widetilde{A}_0$)-($\widetilde{A}_4$) and that
\begin{eqnarray*}\label{th3a1}
  h_0\ {\rm is\ an\ isolated\ zero\ of\ }M
\end{eqnarray*}
with
\begin{eqnarray*}\label{th3b1}
  {\rm ind}\left(h_0,M\right)\not=0.
\end{eqnarray*}
Then, for any $\eps>0$ sufficiently small, system (\ref{ps}) has a
$T$-periodic solution $x_\eps\in C^0([0,T],E)$ and
 \begin{eqnarray*}\label{conve1}
    x_\eps(0)\to\xi_0\quad{\rm as }\ \eps\to 0.
 \end{eqnarray*}
\end{theorem}

\section{Appendix}
\def\theequation{5.\arabic{equation}}\makeatother
\setcounter{equation}{0}

\noindent{\bf Proof of Lemma~\ref{lemimp}.} Let
$\overline{\Phi}_{h,\eps}:E\to E$ be defined by
\begin{equation}\label{vo}
  \overline{\Phi}_{h,\eps}(\xi)=\Phi_{h,\eps}(\pi_h
  \xi)+(I-\pi_h)\xi.
\end{equation}
Observe, that if there exists $r>0,$ $M>0$ and
$\xi:\mathbb{R}^k\times[0,r]\to E$ satisfying
\begin{equation}\label{VOT}
  \xi(\cdot,\eps)\in C^0(V,E),\ \xi(h,\eps)\to\xi(h,0)\ {\rm as}\
  \eps\to 0\mbox{ uniformly in }h\in V,
\end{equation}
such that
\begin{itemize}
\item[b')] $\overline{\Phi}_{h,\eps}(\xi(h,\eps))=0$ for any
$h\in V,$ $\eps\in[0,r],$
\item[c')] $\xi(h,\eps)$ is the only zero of
$\overline{\Phi}_{h,\eps}$ in $B_E(0,r),$
\item[d')] $\|\xi(h,\eps)\|\le M\eps$ for any $h\in V,$
$\eps\in[0,r],$
\end{itemize}
then $\beta(h,\eps)=\pi_h\xi(h,\eps)$ satisfies a), b), c) and d).

\vskip0.2truecm\noindent To prove this assertion from assumption
\ref{cond1} we have
$$
  \overline{\Phi}_{h,0}(0)=\Phi_{h,0}(0)=\widetilde{P}(h,0)=0.
$$
For the derivative $(\overline{\Phi}_{h,0})'(\cdot)$ taking into
account that $\widetilde{P}(h,\cdot)$ acts on $E_h$ we have
$$
  ({\overline{\Phi}}_{h,0})'(0)=\pi_h\widetilde{P}'_\beta(h,0)\pi_h+(I-\pi_h).
$$
Let us show that $(\overline{\Phi}_{h,0})'(0)$ is invertible on $E$
for $h\in V,$ to do this we show that given $b\in E$ there exists a
unique $a_b\in E$ such that
\begin{equation}\label{BE}
(\overline{\Phi}_{h,0})'(0)a_b=b.
\end{equation}
Indeed, applying $I-\pi_h$ to (\ref{BE}) we have
$(I-\pi_h)a_b=(I-\pi_h)b.$ On the other hand, by assumption
\ref{cond2} $\pi_h\widetilde{P}'_\beta(h,0)$ is invertible and thus
applying $\left(\pi_h\widetilde{P}_\beta'(h,0)\right)^{-1}\pi_h$ to
(\ref{BE}) we obtain $\pi_h
a_b=\left(\pi_h\widetilde{P}_\beta'(h,0)\right)^{-1}\pi_h b.$
Therefore the unique solution $a_b$ of (\ref{BE}) is given by
$a_b=\left(\pi_h\widetilde{P}_\beta'(h,0)\right)^{-1}\pi_h
b+(1-\pi_h)b.$ This means that
$((\overline{\Phi}_{h,0})'(0))^{-1}\pi_h$ is continuous in $h.$ Now,
introducing $\overline{P}(h,\xi)=\widetilde{P}(h,\pi_h
\xi)+(I-\pi_h)\xi$ we have that
\begin{itemize}
\item[1')] $\overline{\Phi}_{h,\eps}(\xi)=\overline{P}(h,\xi)+\eps
\widetilde{Q}(h,\xi,\eps),$
\item[2')] $\overline{P}(h,0)=0,$
\item[3')] $\overline{P}'_\xi(h,0)$ is invertible and
$\left(\overline{P}'_\xi(h,0)\right)^{-1}$ is continuous in $h.$
\end{itemize}
Let
$\widehat{\Phi}_{h,\eps}(\xi)=\left(\overline{P}'_{\xi}(h,0)\right)^{-1}\overline{\Phi}_{h,\eps}(\xi).$
Since $\overline{\Phi}_{h,\eps}(\xi)=0$ if and only if
$\widehat{\Phi}_{h,\eps}(\xi)=0$ we aim now at finding a solution
$\xi(h,\eps)$ to $\widehat{\Phi}_{h,\eps}(\xi)=0$ satisfying
properties b'), c') and d'). By assumption \ref{DST1} for any $h\in
\overline{V}$ there exists $r(h)>0$ such that
$$\|I-(\widehat{\Phi}_{\widehat{h},0})'(\xi)\|\le 1/4$$
for any $\|\xi\|\le r(h)$ and any $\widehat{h}\in
B_{\mathbb{R}^k}(h,r(h))\cap\overline{V}.$

\vskip0.2truecm\noindent Since the family $\bigcup_{h\in V}
B_{\mathbb{R}^k}(h,r(h))$ covers the set $\overline{V}$ we can
extract from it a finite subfamily covering $V.$ This implies the
existence of $r>0$ such that
$$\|I-(\widehat{\Phi}_{h,0})'(\xi)\|\le 1/4$$
for any $\|\xi\|\le r$ and any $h\in \overline{V}.$

\vskip0.2truecm\noindent By assumption \ref{DST2} there is $L>0$
such that
$\|(\overline{P}'_{\xi}(h,0))^{-1}(\eps\widetilde{Q}(h,\xi_1,\eps)-\eps\widetilde{Q}(h,\xi_2,\eps))\|\le
\eps L$ for any $h\in V,$ $\xi_1,\xi_2\in B_E(0,1),$ $\eps\in[0,1].$

\vskip0.2truecm\noindent Therefore, $r>0$ can be considered
sufficiently small to have
\begin{equation}\label{US}
  \|\xi_1-\widehat{\Phi}_{h,\eps}(\xi_1)-\xi_2+\widehat{\Phi}_{h,\eps}(\xi_2)\|\le
  (1/2)\|\xi_1-\xi_2\|
\end{equation}
for any $h\in V,$  $\eps\in[0,r],$ $\|\xi_1\|\le r,$ $\|\xi_2\|\le
r$. Therefore, there exists $\xi:V\times[0,r]\to E$ satisfying b')
and c'). It remains to show that $\xi$ satisfies also (\ref{VOT})
and $d').$ Indeed, by using b') and (\ref{US}) for any $h_1,h_2\in
V$ and $\eps\in[0,r]$ we have
  \begin{eqnarray*}
  &&\|\xi(h_1,\eps_1)-\xi(h_2,\eps_2)\|\le\\&\le&
  \|\xi(h_1,\eps_1)-\widehat{\Phi}_{h_2,\eps_2}(\xi(h_1,\eps_1))-\xi(h_2,\eps_2)+
  \widehat{\Phi}_{h_2,\eps_2}(\xi(h_2,\eps_2))\|+ \\
  &&+\|\widehat{\Phi}_{h_2,\eps_2}(\xi(h_1,\eps_1))-\widehat{\Phi}_{h_1,\eps_2}(\xi(h_1,\eps_1))\|+\\&&+
  \|\widehat{\Phi}_{h_1,\eps_2}(\xi(h_1,\eps_1))-\widehat{\Phi}_{h_1,\eps_1}(\xi(h_1,\eps_1))\|\le\\
 & \le &
  (1/2)\|\xi(h_1,\eps_1)-\xi(h_2,\eps_2)\|+\|\widehat{\Phi}_{h_2,\eps_2}(\xi(h_1,\eps_1))-
  \widehat{\Phi}_{h_1,\eps_2}(\xi(h_1,\eps_1))\|+
  \\
  && +|\eps_1-\eps_2|\,\left\|\left(\overline{P}_{\xi}'(h_1,0)\right)^{-1}
  \left(\widetilde{Q}(h_1,\xi(h_1,\eps_1),\eps_2)-
  \widetilde{Q}(h_1,\xi(h_1,\eps_1),\eps_1)\right)\right\|.
  \end{eqnarray*}
Finally the continuity assumptions \ref{DST3}, \ref{DST1} and
\ref{DST2} imply that $\xi$ satisfies (\ref{VOT}) and $d').$ \qed

\end{document}